# BOUNDS ON COVERAGE PROBABILITIES OF THE EMPIRICAL LIKELIHOOD RATIO CONFIDENCE REGIONS[1]

By Min Tsao

*University of Victoria*

This paper studies the least upper bounds on coverage probabilities of the empirical likelihood ratio confidence regions based on estimating equations. The implications of the bounds on empirical likelihood inference are also discussed.

**1. Introduction.** The fact that there is a nontrivial upper bound (less than one) on the coverage probability of an empirical likelihood ratio confidence region is most easily seen through that for the mean. In this case the confidence region is nested within the convex hull of the sample. Thus, regardless of its confidence level, a nontrivial upper bound on its coverage probability is the probability that the convex hull covers the mean.

Several factors affect the value of the upper bound: the underlying distribution, the sample size and the dimension of the mean. In empirical likelihood inference the underlying distribution is not available. Thus, even for the simple case of the mean, the upper bound on coverage probability cannot be determined. Interestingly, however, for a large class of empirical likelihood ratio confidence regions, including those for the mean, the least upper bound on the coverage probability is available. This paper studies this least upper bound and its implications for empirical likelihood inference.

**2. Main results.** To set up notation, consider a parameter of interest $\theta_0$ of a continuous random vector $Y$. Let $Y_1, Y_2, \ldots, Y_n$ be $n$ i.i.d. copies of $Y$. Let $m(Y, \theta) \in \mathcal{R}^k$ be an estimating function for $\theta_0$ that is continuous in $Y$.

Received May 2003; revised October 2003.

[1] Supported by a grant from the Natural Sciences and Engineering Research Council of Canada.

The empirical likelihood ratio function for $\theta_0$ is

$$(2.1) \quad R(\theta) = \sup\left\{\prod_{i=1}^{n} nw_i \;\Big|\; \sum_{i=1}^{n} w_i m(Y_i, \theta) = \underline{0}, w_i \geq 0, \sum_{i=1}^{n} w_i = 1\right\},$$

where $\underline{0}$ is the origin in $\mathcal{R}^k$. See Owen (2001) and Qin and Lawless (1994). The log likelihood ratio $l(\theta)$ is given by $l(\theta) = -2\log R(\theta)$. The empirical likelihood ratio confidence region for $\theta_0$ is given by

$$(2.2) \quad \mathcal{C}_r = \{\theta | l(\theta) < r\},$$

where $r$ is a finite quantity determined by the desired confidence level through the method of calibration of choice. Throughout this paper the sample size $n$ and the dimension of the estimating function $k$ are assumed fixed unless we specify them to be otherwise.

Denote by $\mathcal{H}(m(Y_1, \theta_0), m(Y_2, \theta_0), \ldots, m(Y_n, \theta_0))$ the convex hull of $m(Y_i, \theta_0)$. Because $l(\theta_0)$ is finite if and only if $\underline{0}$ is in the interior of the convex hull, event $\{\theta_0 \in \mathcal{C}_r\}$ implies $\{\underline{0} \in \mathcal{H}(m(Y_1, \theta_0), m(Y_2, \theta_0), \ldots, m(Y_n, \theta_0))\}$. Thus,

$$(2.3) \quad P(\theta_0 \in \mathcal{C}_r) < P[\underline{0} \in \mathcal{H}(m(Y_1, \theta_0), m(Y_2, \theta_0), \ldots, m(Y_n, \theta_0))].$$

Further, $P(\theta_0 \in \mathcal{C}_r)$ is a monotone increasing function of $r$ and

$$(2.4) \quad \lim_{r \to +\infty} P(\theta_0 \in \mathcal{C}_r) = P[\underline{0} \in \mathcal{H}(m(Y_1, \theta_0), m(Y_2, \theta_0), \ldots, m(Y_n, \theta_0))].$$

Hence, the bound in the right-hand side of (2.3) is the least upper bound on the coverage probability of the confidence region (2.2) associated with the particular $m(Y, \theta_0)$. This bound, however, is in general not available because the distribution of $m(Y, \theta_0)$ is not available. We consider instead the least upper bound $B$,

$$B = \sup\{P(\theta_0 \in \mathcal{C}_r)\}$$
$$= \sup\{P[\underline{0} \in \mathcal{H}(m(Y_1, \theta_0), m(Y_2, \theta_0), \ldots, m(Y_n, \theta_0))]\},$$

where the supremum is taken over all empirical likelihood ratio confidence regions based on estimating equations (2.1) and (2.2), or equivalently, all meaningful $m(Y, \theta_0)$ and $r$.

In order to find $B$ without having to characterize the set of all meaningful $m(Y, \theta_0)$, let $X_1, X_2, \ldots, X_n$ be i.i.d. copies of an arbitrary continuous random vector $X$ in $\mathcal{R}^k$ and denote by $\mathcal{H}(X_1, X_2, \ldots, X_n)$ their convex hull. Consider $b(k, n)$ given by

$$(2.5) \quad b(k, n) = \sup_{X}\{P[\underline{0} \in \mathcal{H}(X_1, X_2, \ldots, X_n)]\},$$

where the supremum is taken over all possible continuous random vectors in $\mathcal{R}^k$. We claim that (i) $b(k, n)$ is attained at an $X$ if and only if the distribution of its projection on the unit sphere $X^p$ is symmetric with respect to $\underline{0}$ and (ii) $b(k, n)$ is the least upper bound $B$. Once (i) is established, (ii) then follows from the fact that any $m(Y, \theta_0)$ is a special case of $X$ and that $b(k, n)$ is attained at a special $m(Y, \theta_0)$ for, say, the empirical likelihood inference for the mean of the uniform distribution on the unit sphere in $\mathcal{R}^k$. To see the latter point, since $Y$ is uniform on the unit sphere, $\theta_0 = E(Y) = \underline{0}$ and $m(Y, \theta_0) = Y - \theta_0 = Y$. Hence, this $m(Y, \theta_0)$ and its projection are both symmetric with respect to $\underline{0}$. To prove claim (i), we need the following lemma.

LEMMA 1. *For any continuous $X$ in $\mathcal{R}^k$, let $v_i = \|X_i\|_2$ and without loss of generality assume $v_i > 0$. Let $X_i^p = v_i^{-1} X_i$ be the projection of $X_i$ on the unit sphere. Then*

$$P\{\underline{0} \in \mathcal{H}(X_1, X_2, \ldots, X_n)\} = P\{\underline{0} \in \mathcal{H}(X_1^p, X_2^p, \ldots, X_n^p)\}.$$

PROOF. It suffices to show that $\underline{0} \notin \mathcal{H}(X_1, X_2, \ldots, X_n)$ if and only if $\underline{0} \notin \mathcal{H}(X_1^p, X_2^p, \ldots, X_n^p)$. The convex hull $\mathcal{H}(X_1, X_2, \ldots, X_n)$ does not contain $\underline{0}$ if and only if all $X_i$ are on one side of a hyperplane through $\underline{0}$. All $X_i$ are on one side of a hyperplane through $\underline{0}$ if and only if their projections $X_i^p$ are on one side of a hyperplane through $\underline{0}$. All $X_i^p$ are on one side of a hyperplane if and only if their convex hull $\mathcal{H}(X_1^p, X_2^p, \ldots, X_n^p)$ does not contain $\underline{0}$. Thus the lemma. $\square$

Claim (i) implies that $b(k, n) = P\{\underline{0} \in \mathcal{H}(U_1, U_2, \ldots, U_n)\}$, where $U_1, U_2, \ldots, U_n$ are i.i.d. copies of a uniform random vector $U$ supported on the unit sphere in $\mathcal{R}^k$. We now prove this claim for $k = 1, 2$.

THEOREM 1. *Let $k = 1, 2$ and $n > k$. For any continuous $X$ in $\mathcal{R}^k$, we have*

$$(2.6) \qquad P\{\underline{0} \in \mathcal{H}(X_1, X_2, \ldots, X_n)\} \leq P\{\underline{0} \in \mathcal{H}(U_1, U_2, \ldots, U_n)\}.$$

*Further, equality holds if and only if the distribution of the projection of $X$ on the unit sphere $X^p$ is symmetric with respect to $\underline{0}$.*

PROOF. By Lemma 1 we only need to show that (2.6) holds for all continuous $X$ supported on the unit sphere. Thus, we assume without loss of generality that $X$ is supported on the unit sphere. Under this assumption, the symmetry condition on $X^p$ in Theorem 1 is equivalent to the symmetry condition on $X$ itself.

For $k = 1$, the unit sphere and, thus, the support of $X$ degenerates into $\{-1, 1\}$. Let $p = P\{X = 1\}$. Then

$$P\{\underline{0} \in \mathcal{H}(X_1, X_2, \ldots, X_n)\} = 1 - p^n - (1 - p)^n.$$

Theorem 1 amounts to the simple observation that function $1 - p^n - (1-p)^n$ attains its unique maximum at $p = 1/2$ which corresponds to the uniform distribution on $\{-1, 1\}$, the only symmetric distribution on $\{-1, 1\}$.

For $k = 2$, let $X$ be a continuous random variable on the unit circle $(0 \leq X < 2\pi)$ and for simplicity assume that its density $f(x)$ is continuous on the circle. Define

$$G(x) = \int_x^{x+\pi} f(y) \, dy,$$

where $f(x) = f(2\pi + x)$. For $X_1, \ldots, X_{j-1}, X_{j+1}, \ldots, X_n$, denote the event that they are in the half-circle $(X_j, X_j + \pi)$ by $A_j$. If $X_j > \pi$, this half-circle represents the union of $(X_j, 2\pi)$ and $[0, X_j - \pi)$. Since $X_i$ are i.i.d., we have for $j = 1, 2, \ldots, n$

$$P\{A_j\} = \int_0^{2\pi} f(x)[G(x)]^{n-1} \, dx.$$

Further, $A_i \cap A_j = \phi$ for $i \neq j$, where $\phi$ denotes the empty set, and

$$\{\underline{0} \notin \mathcal{H}(X_1, X_2, \ldots, X_n)\} = \bigcup_{i=1}^n A_i.$$

It follows that for any $n \geq 1$,

$$
\begin{aligned}
(2.7) \quad P\{\underline{0} \notin \mathcal{H}(X_1, X_2, \ldots, X_n)\} &= \sum_{i=1}^n P\{A_i\} \\
&= n \int_0^{2\pi} f(x)[G(x)]^{n-1} \, dx.
\end{aligned}
$$

Noting that $P\{A_j\}$ equals the probability that $X_1, \ldots, X_{j-1}, X_{j+1}, \ldots, X_n$ are in the half-circle $(X_j - \pi, X_j)$, an equivalent expression for $P\{\underline{0} \notin \mathcal{H}(X_1, X_2, \ldots, X_n)\}$ is

$$
\begin{aligned}
(2.8) \quad P\{\underline{0} \notin \mathcal{H}(X_1, X_2, \ldots, X_n)\} &= \sum_{i=1}^n P\{A_i\} \\
&= n \int_0^{2\pi} f(x)[G(x - \pi)]^{n-1} \, dx.
\end{aligned}
$$

Adding up (2.7) and (2.8) gives another expression for $P\{\underline{0} \notin \mathcal{H}(X_1, X_2, \ldots, X_n)\}$,

$$
\begin{aligned}
(2.9) \quad & P\{\underline{0} \notin \mathcal{H}(X_1, X_2, \ldots, X_n)\} \\
&= \frac{n}{2} \int_0^{2\pi} f(x)\{[G(x)]^{n-1} + [G(x - \pi)]^{n-1}\} \, dx.
\end{aligned}
$$

To see that the equality in (2.6) holds if the distribution of $X$ is symmetric with respect to $\underline{0}$, note that the distribution is symmetric if and only if

$G(x) = 1/2$ for all $x \in [0, 2\pi)$. This and (2.7) imply that for all symmetric $X$, including $U$,

$$(2.10) \qquad P\{\underline{0} \in \mathcal{H}(X_1, X_2, \ldots, X_n)\} = 1 - n(1/2)^{n-1}.$$

To show that the inequality in (2.6) holds strictly if the distribution of $X$ is not symmetric and, thus, it also must be symmetric if the equality holds, first note that for any $n \geq 1$ and $p \in [0, 1]$, the function $h(p) = p^{n-1} + (1-p)^{n-1}$ achieves its unique minimum at $p = 1/2$ and this minimum is $h(1/2) = (1/2)^{n-2}$. Since $G(x), G(x - \pi) \geq 0$ and $G(x) + G(x - \pi) = 1$, for any $n \geq 1$,

$$(2.11) \qquad (1/2)^{n-2} \leq [G(x)]^{n-1} + [G(x - \pi)]^{n-1}.$$

If the distribution of $X$ is not symmetric, $G(x)$ cannot be $1/2$ for all $x \in [0, 2\pi)$. Further, $G(x)$ is continuously differentiable. There exists an open subinterval of $[0, 2\pi)$ in which $G(x) \neq 1/2$ and $G'(x) < 0$. Over this subinterval $f(x) > 0$ and the inequality in (2.11) holds strictly. Multiply both sides of (2.11) by $f(x)$ and then integrate from 0 to $2\pi$. We have

$$(2.12) \qquad (\tfrac{1}{2})^{n-1} < \tfrac{1}{2} \int_0^{2\pi} f(x)\{[G(x)]^{n-1} + [G(x - \pi)]^{n-1}\}\, dx,$$

where the left-hand side is strictly smaller than the right-hand side because of the subinterval. It follows from (2.9), (2.12) and (2.10) that for an $X$ that is not symmetric,

$$P\{\underline{0} \in \mathcal{H}(X_1, X_2, \ldots, X_n)\} < 1 - n(1/2)^{n-1}$$
$$= P\{\underline{0} \in \mathcal{H}(U_1, U_2, \ldots, U_n)\}. \qquad \square$$

For $k \geq 3$, a proof of (2.6) has eluded us so far due to difficulties in finding an analytic expression for $P\{\underline{0} \notin \mathcal{H}(X_1, X_2, \ldots, X_n)\}$ for a general $X$ in high dimensions. Thus, claim (i) has been proved for only $k \leq 2$. We conjecture that claim (i) holds for all $k$. The rest of our discussion assumes this conjecture holds so that $b(k, n) = P\{\underline{0} \in \mathcal{H}(U_1, U_2, \ldots, U_n)\}$ for all $k$. Wendel (1962) gives a formula for $P\{\underline{0} \notin \mathcal{H}(U_1, U_2, \ldots, U_n)\}$ which leads to the following expression for $b(k, n) = P\{\underline{0} \in \mathcal{H}(U_1, U_2, \ldots, U_n)\}$: for any $n > k$,

$$(2.13) \quad b(k, n) = 1 - \left\{ \binom{n-1}{0} + \binom{n-1}{1} + \cdots + \binom{n-1}{k-1} \right\} 2^{-(n-1)}.$$

It is interesting to note that, by (2.13), when the sample size is twice as much as the dimension, the value of the least upper bound $b(k, 2k)$ equals 0.5. Theorem 2 further explores the implications of (2.13).

THEOREM 2. *Denote by $[x]$ the largest integer smaller than $x$. For any $n > k$,*

(a) $b(k, n + 1) > b(k, n)$ and $b(k, n) > b(k + 1, n),$ and
(b) for any $\varepsilon \in (0, 0.5),$ $b([\varepsilon n], n) \to 1$ and $b([(1 - \varepsilon)n], n) \to 0$ as $n \to \infty.$

PROOF. The inequalities in (a) follow easily from (2.13). To see (b) is true, consider the binomial random variable $X \sim \text{Bin}(1/2, n - 1)$. Denote by $Z$ the standard normal random variable. By (2.13) we have

$$b([\varepsilon n], n) = 1 - P\{X \leq [\varepsilon n] - 1\}$$
$$\sim 1 - P\left\{Z \leq \frac{[\varepsilon n] - 1 - (n - 1)/2}{\sqrt{n - 1}/2}\right\}.$$

The right-hand side and, thus, $b([\varepsilon n], n)$ go to one when $n$ goes to infinity. Similarly, $b([(1 - \varepsilon)n], n)$ goes to zero as $n$ goes to infinity. $\square$

3. **Concluding remarks.** The least upper bound $b(k, n)$ may be surprisingly small when the ratio $n/k$ is small. Table 1 shows values of the bound at various combinations of $k$ and $n$. When this ratio is small and an empirical likelihood ratio confidence region of a high confidence level is desired, it is essential that the bound be computed to see if such a high confidence level is impossible. We have come across examples in the literature where regions of impossibly high confidence levels were computed. Practitioners need to be aware of the bound.

For any fixed $n$, the bound $b(k, n)$ is a strictly decreasing function of $k$. When the sample size $n$ is not large, practitioners need to be aware of the negative impact of incorporating extra information about the parameter that will increase the dimension of the estimating equation $k$: (i) high confidence levels may become unachievable and (ii) continuous approximations to the finite sample distribution of the empirical log likelihood ratio may also become less accurate. The latter may diminish the benefit of incorporating the extra information and may, for some cases where $n$ is not large, result in a loss in coverage accuracy for the empirical likelihood ratio confidence region [Tsao (2004)].

The method of empirical likelihood has been applied to some very high-dimensional problems and there is increasing interest in the asymptotic behavior of the empirical log likelihood ratio when the sample size $n$ and the dimension of the estimating equation $k$ both tend to infinity. By Theorem 2, when $n \leq \gamma k$ for some constant $\gamma \in (1, 2)$ and $n$ goes to infinity, the distribution of the empirical log likelihood ratio $l(\theta_0)$ will degenerate into a point mass at infinity. There are no meaningful confidence regions of the form (2.2) in this case.

TABLE 1
*Bounds for some combinations of $n$ and $k$, $r = n/k$*

| k | r = 2 | r = 3 | r = 4 | r = 5 | r = 6 | r = 7 | r = 8 |
|---|-------|-------|-------|-------|-------|-------|-------|
| 1 | 0.5000 | 0.7500 | 0.8750 | 0.9375 | 0.9688 | 0.9844 | 0.9922 |
| 2 | 0.5000 | 0.8125 | 0.9375 | 0.9805 | 0.9941 | 0.9983 | 0.9995 |
| 5 | 0.5000 | 0.9102 | 0.9904 | 0.9992 | 0.9999 | 1.0000 | 1.0000 |

On related future research problems, we note that in light of the lack of awareness of the bounds, a method of calibration which automatically respects the bounds may be helpful. Tsao (2004) contains some preliminary results on one such method. It may be possible to derive similar bounds for certain classes of empirical likelihood ratio confidence regions outside of the estimating equation framework (2.1) and (2.2). The conjecture that claim (i) holds for all $k$ is another interesting question that we are still working on.

To conclude, while trying to determine the value of the bound, based on derivations for $k = 1, 2$ and some asymptotic observations we had communicated to several colleagues the conjecture that for any $k$ and $n > k$,

$$P\{\underline{0} \notin \mathcal{H}(U_1, U_2, \ldots, U_n)\}$$
$$= \left\{ \binom{n-1}{0} + \binom{n-1}{1} + \cdots + \binom{n-1}{k-1} \right\} 2^{-(n-1)}.$$

We are indebted to Professor Qi-Man Shao who brought to our attention related work by J. G. Wendel, B. Efron and others. Efron (1965) appears to be the first to give formulae (2.7) and (2.8). Wendel (1962) has already noted and proved the conjecture. Citing connections to L. J. Savage, R. E. Machol and D. A. Darling, Wendel (1962) also gives an interesting historical note on the origin of the conjecture.

**Acknowledgments.** The author would like to thank Professors Art B. Owen, Qi-Man Shao, Song Xi Chen and Zhidong Bai for helpful discussions and two referees for helpful comments. A part of this work was done while the author was visiting the Department of Statistics and Applied Probability at the National University of Singapore.

Department of Mathematics
and Statistics
University of Victoria
Victoria, British Columbia
Canada V8W 3P4
e-mail: tsao@math.uvic.ca